\documentclass[conference]{IEEEtran}
\IEEEoverridecommandlockouts
\usepackage{cite}
\usepackage{amsmath,amssymb,amsfonts}
\usepackage{algorithmic}
\usepackage{graphicx}
\usepackage{textcomp}
\usepackage{xcolor}
\def\BibTeX{{\rm B\kern-.05em{\sc i\kern-.025em b}\kern-.08em
    T\kern-.1667em\lower.7ex\hbox{E}\kern-.125emX}}
\begin{document}

\font\myfont=cmr18 at 20pt
\title{{\myfont  A modified limited memory Nesterov's accelerated quasi-Newton}
  }

\author{\IEEEauthorblockN{S. Indrapriyadarsini\IEEEauthorrefmark{2},
Shahrzad Mahboubi\IEEEauthorrefmark{3}, Hiroshi Ninomiya\IEEEauthorrefmark{3},
Takeshi Kamio\IEEEauthorrefmark{4} and
Hideki Asai\IEEEauthorrefmark{2}}
\IEEEauthorblockA{Shizuoka University \IEEEauthorrefmark{2}, 
    Shonan Institute of Technology\IEEEauthorrefmark{3}, 
    Hiroshima City University\IEEEauthorrefmark{4}\\
\IEEEauthorrefmark{2}indra.ipd@gmail.com}}

\maketitle
\section{Introduction}
The Nesterov's accelerated quasi-Newton (L)NAQ \cite{lnaq} method has shown to accelerate the conventional (L)BFGS quasi-Newton method using the Nesterov's accelerated gradient in several neural network (NN) applications. However, the calculation of two gradients per iteration increases the computational cost. An approximation to the Nesterov's accelerated gradient was proposed in \cite{pricai}. This paper extends the study in \cite{pricai} by applying the approximation to limited memory NAQ. 
\vspace{-1mm}
\section{Proposed Algorithm (L-MoQ)}

\begin{figure}[!b]
\begin{center}
\includegraphics[width=8cm]{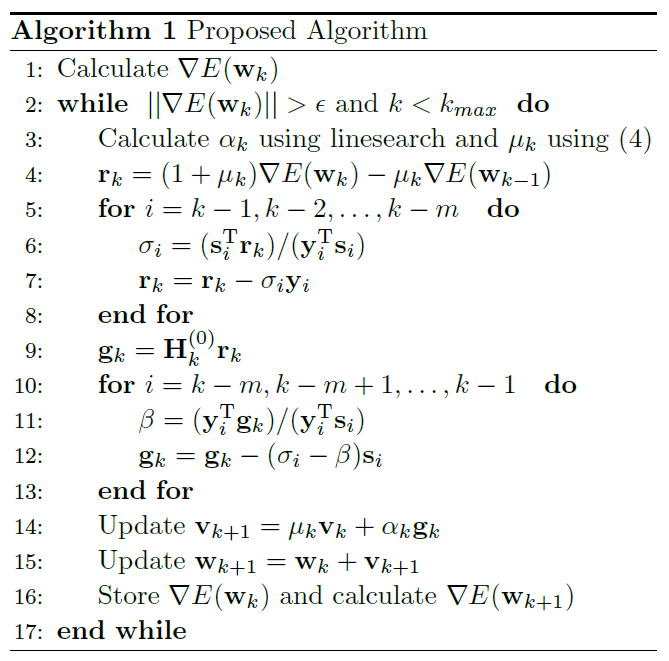}
\end{center}
\vspace{-1.5mm}
\end{figure}
\vspace{-1mm}

The weight update in NAQ is given as ${\bf w}_{k+1} = {\bf w}_{k} + \mu_k{\bf v}_{k} + \alpha_k {\bf g}_{k},$ where ${\bf g}_{k}$ is the search direction as  in (1).  
\vspace{-3mm}
\begin{equation}
{\bf g}_{k} =  {\bf H}_{k}  {\bf \nabla E(w}_{k}+\mu_k{\bf v}_{k}).
\vspace{-1mm}
\end{equation}
MoQ [2] showed that (1) can be approximated as 
\vspace{-1mm}
\begin{equation}
{\boldsymbol { \bf g}}_k=  {\bf H}_{k} [ (1+ \mu_k) \nabla {\bf E ( w}_k) - \mu_k \nabla {\bf  E(w}_{k-1}) ].
\vspace{-1.5mm}
\end{equation}
The $\mu_k$ value is updated using (3) where $\theta_{k+1}$ is  obtained by solving (4) with $\theta_{0} = 1$ and $\gamma = 10^{-5}$.
\vspace{-2mm}
\begin{equation}
\mu_k = \theta_k(1-\theta_k)/({\theta_k}^2 + \theta_{k+1}),
\vspace{-1.5mm}
\end{equation}
\vspace{-4mm}
\begin{equation}
{\theta_{k+1}}^2 = (1+\theta_{k+1}){\theta_k}^2+\gamma \theta_{k+1}.
\vspace{-1.5mm}
\end{equation}

 In this paper, we propose limited-memory MoQ (L-MoQ) in which (2) is evaluated by the two-loop recursion using the last $m$ curvature pairs $({\bf s} , {\bf y})$ given by, 
 \vspace{-4mm}
\begin{equation*}
 ~~~~~~~~{\bf s}_k = {\bf w}_{k+1} - ({\bf w}_k + \mu_k {\bf v}_{k}) ~~~{\rm and }~~
 \vspace{-1.5mm}
\end{equation*}
\begin{equation}
  {\bf y}_k = \nabla {\bf E(w}_{k+1}) - [(1+ \mu_k) \nabla {\bf E(w}_k) - \mu_k \nabla {\bf E(w}_{k-1})] ).
  \vspace{-1.5mm}
\end{equation}

\vspace{-3.2mm}
\section{Simulation Results}
\vspace{-3.2mm}
\begin{equation*}
 f(x_1 \ldots x_n) = \frac{\pi}{n} \Bigr \{ \sum_{i=1}^{n-1} [(x_i-1)^2(1+10~ {\rm sin}^2(\pi x_{i+1}))]~~~~
\end{equation*}
\vspace{-3.2mm}
\begin{equation}
  ~~~~~~~~~~ + 10~ {\rm sin}^2(\pi x_1) +(x_n-1)^2 \Bigl \}, x_i \in [-4,4], \forall i.
\end{equation}

The performance of L-MoQ is evaluated on Levy function (6) using a $5-50-1$ NN with $k_{\rm max}=10000$, $\epsilon=10^{-6}$ and $m=16$. The number of parameters is $d= 351$. Fig. 1 and Table 1 show the average results of 50 trials. The number of function and gradient evaluations are denoted as $fev$ and $gev$, respectively. The results confirm that the proposed L-MoQ is a good approximation to L-NAQ while having fewer $gev$ and maintaining the same computational cost of L-BFGS i.e., $nd+4md+2d+\zeta nd$ and storage cost of $(2m+1)d$. 


\begin{table}[!t]
\begin{center}
\caption{Summary of results averaged over 50 trials.}\label{table}\vspace{-2mm}
\begin{tabular}{|c|c|c|c|c|c|}
\hline 
Method 	& $E{\bf (w)}  $ &	iters  ($k$)& $fev$	&$gev$ & time(s)\\
\hline	
L-BFGS & 0.000091   & 10000	   & 28398	&  10001   & 81.95 \\ \hline
L-NAQ 	 & 0.000025    &	9927		& 20848	   &  19854  & 95.21 \\ \hline
{\bf L-MoQ }	 & \bf{0.000022}	  & {\bf 9961}  & {\bf 20918}		& {\bf 9962}	   & {\bf 73.73} \\

\hline
\end{tabular}
\end{center}
\vspace{-3mm}
\end{table}
\vspace{-3mm}
\begin{figure}[!htb]
\begin{center}
\includegraphics[width=8.2cm]{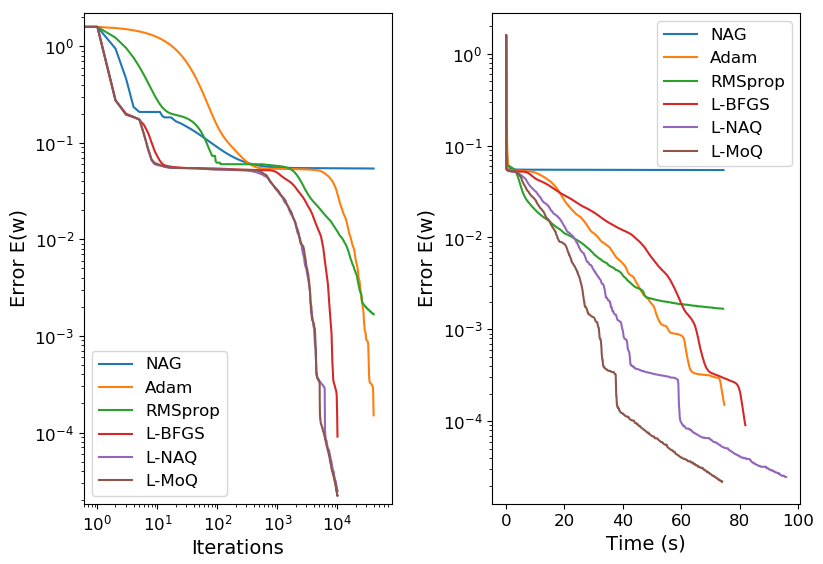}
\end{center}
\vspace{-1.5mm}
\caption{Average training error over 50 trials. } 
\end{figure}

\end{document}